\input amstex
\documentstyle{amsppt}
 
\def\binrel@#1{\setbox\z@\hbox{\thinmuskip0mu
\medmuskip\m@ne mu\thickmuskip\@ne mu$#1\m@th$}%
\setbox\@ne\hbox{\thinmuskip0mu\medmuskip\m@ne mu\thickmuskip
\@ne mu${}#1{}\m@th$}%
\setbox\tw@\hbox{\hskip\wd\@ne\hskip-\wd\z@}}
\def\overset#1\to#2{\binrel@{#2}\ifdim\wd\tw@<\z@
\mathbin{\mathop{\kern\z@#2}\limits^{#1}}\else\ifdim\wd\tw@>\z@
\mathrel{\mathop{\kern\z@#2}\limits^{#1}}\else
{\mathop{\kern\z@#2}\limits^{#1}}{}\fi\fi}
\def\underset#1\to#2{\binrel@{#2}\ifdim\wd\tw@<\z@
\mathbin{\mathop{\kern\z@#2}\limits_{#1}}\else\ifdim\wd\tw@>\z@
\mathrel{\mathop{\kern\z@#2}\limits_{#1}}\else
{\mathop{\kern\z@#2}\limits_{#1}}{}\fi\fi}
\def\circle#1{\leavevmode\setbox0=\hbox{h}\dimen@=\ht0
\advance\dimen@ by-1ex\rlap{\raise1.5\dimen@\hbox{\char'27}}#1}
\def\sqr#1#2{{\vcenter{\hrule height.#2pt
     \hbox{\vrule width.#2pt height#1pt \kern#1pt
       \vrule width.#2pt}
     \hrule height.#2pt}}}

\def\force{\hbox{$\|\hskip-2pt\hbox{--}$\hskip2pt}}  
 
\baselineskip 24pt

\define\pmf{\par\medpagebreak\flushpar}
\define\pbf{\par\bigpagebreak\flushpar}
\define\k{\kappa}
\define\a{\alpha}

\topmatter
\title
Full reflection of stationary sets at regular cardinals
\endtitle
\author
Thomas Jech and Saharon Shelah
\endauthor
\thanks 
The first author was supported by NSF and by a Fulbright grant.
He wishes to express his gratitude to the Hebrew University
for its hospitality.\endgraf
The second author wishes to thank the
BSF. This paper is No. 383 on his publication list\endthanks
\affil The Pennsylvania State University\\The Hebrew University \endaffil
\address \endgraf Department of Mathematics, The Pennsylvania State University,
University Park, PA 16803, USA \endgraf Institute of Mathematics,
The Hebrew University, Jerusalem, Israel \endaddress
\endtopmatter
 
\document
\subhead
0. Introduction
\endsubhead
 
A stationary subset $S$ of a regular uncountable cardinal $ \k$
{\it reflects
fully at regular cardinals} if for every stationary set $T
\subseteq \k$ of higher order consisting of regular cardinals
there exists an $ \a \in T$ such that $S \cap \a $ is a stationary
subset of $\a$.
We prove that the Axiom of Full Reflection which states that
every stationary set reflects
fully at regular cardinals, together with the
existence of $n$-Mahlo cardinals
is equiconsistent with the existence of $ \Pi^1_n$-indescribable
cardinals.  We
also state the appropriate generalization for greatly Mahlo
cardinals.

\document
\subhead
1. Results 
\endsubhead
 
It has been proved [7], [3] that reflection of stationary sets is
a large
cardinal property.  We address the question of what is the
largest possible
amount of reflection.  Due to complications that arise at
singular ordinals, we deal in this paper exclusively with
reflection at regular cardinals.  (And so we deal
with stationary subsets of cardinals that are at least Mahlo
cardinals.)\footnote{If
$\k \ge \aleph_3$ then there exist stationary sets $S \subseteq
\{ \alpha < \k : \text{cf }
\alpha = \aleph_0 \}$ and $T \subseteq \{ \beta < \k : \text{cf } \beta =
\aleph_1 \}$, such that $S$ does not
reflect at any $\beta \in T$.}
 
If $S$ is a stationary subset of a regular uncountable cardinal
$\k,$ then the {\it trace of $S$} is the set
$$ Tr (S) = \{ \a < \k : \; \; S \cap \a \; \text{is stationary
in}  \; \a \} $$
(and we say that $S$ {\it reflects at} $\a$).  If $S$ and $T$ are
both stationary,
we define
$$ S < T \;\; \text{if for almost all} \; \a \in T, \; \; \a \in
Tr(S) $$
and say that $S$ {\it reflects fully} in $T$.
(Throughout the paper, ``for almost all" means ``except for a
nonstationary set of
points").  As proved in [4],  $< $  is a well founded relation;
the {\it order}
$o(S) $ of a stationary set is the rank of $S$ in this relation.
 
If the trace of $S$ is stationary, then clearly $o(S) < o
(Tr(S)).$  We say that $S$
{\it reflects fully at regular cardinals} if its trace meets
every stationary set $T$ of regular cardinals such that $
o(S) < o(T).$  In other words, if for all stationary sets $T$ of
regular cardinals,
$$ o(S) < o(T) \; \text{implies} \; S < T. $$
 
\proclaim
{Axiom of Full Reflection for $\k$}  Every stationary subset of $\k$ reflects
fully at regular cardinals.
\endproclaim
\pmf
 
In this  paper we investigate full reflection together with the
existence of cardinals in the Mahlo hierarchy.  Let $Reg$ be the
set of all regular
limit cardinals $ \a < \k,$ and for each $\eta < \k^+ $ let
$$ E_{ \eta} = Tr^{ \eta} (Reg) - Tr^{ \eta + 1} (Reg) $$
(cf. [2]), and call $ \k \; \;  \eta$-{\it Mahlo} where $ \eta \le
\k^+$ is the least
$ \eta$ such that $E_{ \eta}$ is nonstationary.  In particular,
\flushpar
$E_0 = $ inaccessible non Mahlo cardinals
\flushpar
$E_1 = $ 1-Mahlo cardinals, etc.
\flushpar
We also denote
\flushpar
$E_{-1} = Sing =$ the set of all singular ordinals $\alpha < \k$.
\flushpar
It is well known [4] that each $E_{\eta}, $ the $ \eta$th
{\it canonical}
stationary set is equal (up to the equivalence almost everywhere)
to the set
$$ \{ \a < \k : \; \; \a \; \text{ is} \; f_{ \eta} ( \a)\text{-Mahlo}
\} $$
where $f_{ \eta} $ is the {\it canonical} $ \eta$th function.  A $
\k^+${\it -Mahlo}
cardinal $ \k$ is called {\it greatly Mahlo} [2].
 
If $ \k$ is less than greatly Mahlo (or if it is greatly Mahlo
and the canonical
stationary sets form a maximal antichain) then Full Reflection
for $ \k$ is equivalent to the statement
\pmf
 
{\sl For every $\eta \ge  -1$, every stationary $S \subseteq E_{
\eta}$ reflects almost
everywhere in $E_{\eta + 1}.$}
\pmf
(Because then the trace of $S$ contains almost
all of every $E_{\nu},$ $ \nu > \eta).$
 
The simplest case of full reflection is when $ \k$ is 1-Mahlo;
then full
reflection  states that every stationary $S \subseteq Sing$
reflects at almost every
$ \a \in E_0.$  We will show that this is equiconsistent
with the existence of a weakly compact cardinal.  More generally,
we shall prove that full reflection
together with the existence of $n$-Mahlo cardinals is
equiconsistent with the
existence of $\Pi^1_n-$indescribable cardinals.
 
To state the general theorem for cardinals higher up in the Mahlo
hierarchy, we first
give some definitions.  We assume that the reader is familiar
with $\Pi^1_n$-indescribability.
A ``formula" means a formula of second order logic for $\langle
V_\k, \in \rangle$.
\definition{Definition}  (a) A formula is $\Pi^1_{\eta +1}$ if it is of
the form $\forall
X \neg \varphi$ where $\varphi$ is a $\Pi^1_\eta$ formula.
 
(b)  If $\eta < \k^+$ is a limit ordinal, a formula is
$\Pi^1_\eta$ if it is of the form
$\exists \nu < \eta \,\varphi (\nu, \cdot)$ where $\varphi (\nu, \cdot)$
is a $\Pi^1_\nu$ formula.
\enddefinition
 
For $\alpha \le \k$ and $\eta < \k^+$ we define the satisfaction
relation
$\langle V_\alpha,\in\rangle \models \varphi$ for $\Pi^1_\eta$
formulas in the
obvious way, the only difficulty arising for limit $\eta$, which
is handled as follows:
$$ \langle V_\alpha,\in\rangle \models \exists \nu < \eta \,\, \varphi
(\nu, \cdot) \;\; \text{ if } \;\;
\exists \nu < f_\eta (\alpha) \langle V_\alpha,\in\rangle \models
\varphi (\nu, \cdot) $$
where $f_\eta$ is the $\eta$th canonical function.
 
\definition{Definition}  $\k$ is $\Pi^1_\eta$-{\it indescribable} ($\eta
< \k^+$) if for
every $\Pi^1_\eta$ formula $\varphi$ and every $Y \subseteq V_\k$,
if $\langle V_\k,\in\rangle \models \varphi (Y)$ then
there exists some $\alpha < \k$ such that $\langle V_\alpha, E
\rangle \models
\varphi (Y \cap V_\alpha)$.
 
$\k$ is $\Pi^1_{\k^+}$-{\it indescribable} if it is $\Pi^1_\eta$-indescribable
for all $\eta < \k^+$.
\enddefinition
\pmf
\proclaim{Theorem A}  Assuming the Axiom of Full Reflection for $\k,$
we have for every $\eta \le (\k^+)^L$:
Every $\eta$-Mahlo cardinal is $\Pi^1_\eta$-indescribable in
$L$.
\endproclaim
\proclaim{Theorem B}  Assume that the ground model satisfies $ V = L.$
There is a generic extension
$V[G]$ that preserves cardinals and cofinalities (and satisfies
GCH)
such that for every cardinal $ \k $ in $V$ and every $\eta \le
\k^+$:
\roster
\item{(a)} If $ \k$ is $\Pi^1_\eta$-indescribable in $V$
then $\k$ is $\eta$-Mahlo in $V[G].$
\item{(b)} $V[G]$ satisfies the Axiom of Full Reflection.
\endroster
\endproclaim
\pbf
 
\subhead{2. Proof of Theorem A} \endsubhead
 
Throughout this section we assume full reflection.  The theorem
is proved by
induction on $\k$.  We shall give the proof for the finite case
of the Mahlo
hierarchy; the general case requires only minor modifications.
 
Let $F^\k_0$ denote the club filter on $\k$ in $L$, and for $n > 0$,
let
$F^\k_n$ denote the $\Pi^1_n$ filter on $\k$ in $L$, i.e. the
filter on
$P (\k) \cap L$ generated by the sets $\{ \alpha < \k: L_\alpha
\models \varphi \}$
where $\varphi$ is a $\Pi^1_n$ formula true in $L_\kappa$.  If $\k$ is 
$\Pi^1_n$-indescribable then $F^\k_n$ is a proper filter.  The
$\Pi^1_n$ ideal on $\k$ is the dual of $F^\k_n$.
 
By induction on $n$ we prove the following lemma which implies
the theorem.
 
\proclaim{Lemma 2.1}  Let $A \in L$ be a subset of $\k$ that is in the
$\Pi^1_n$ ideal.  Then
$A \cap E_{n-1}$ is nonstationary.
\endproclaim
 
To see that the Lemma implies Theorem A, let $n \ge 1$, and
letting $A = \k$, we have the implication
$$ \k \; \text{is in the} \; \Pi^1_n \; \text{ideal in} \;
L \Rightarrow E_{n-1} \; \text{is
nonstationary}, $$
and so
$$ \k \; \text{is not} \; \Pi^1_n \text{-indescribable in} \; L
\Rightarrow \k \; \text{is not}\; n\text{-Mahlo}.$$
 
\demo{Proof}  The case $n = 0$ is trivial (if $A$ is
nonstationary in $L$ then $A \cap Sing$ is nonstationary).  Thus
assume that the statement is true for $n$, for all $\lambda \le
\k$, and
let us prove it for $n +1$ for $\k$.  Let $A$ be a subset of $\k,$
$A \in L$, and let $\varphi$ be a $\Pi^1_n$ formula
such that for all $\alpha \in A$ there is some $X_\alpha \in L,$
$X_\alpha \subseteq \alpha,$ such
that $L_\alpha \models \varphi (X_\alpha)$.  Assuming that
$A \cap E_n$ is stationary, we
shall find an $X \in L,$ $X \subseteq \k$, such that $L_\k \models
\varphi (X)$.  Let
$B \supseteq A$ be the set
$$ B = \{ \alpha < \k : \exists X \in L \,\, L_\alpha \models
\varphi (X) \}, $$
and for each $\alpha \in B$ let $X_\alpha$ be the least such $X$
(in $L$).  For each
$\alpha \in B,$ $X_\alpha \in L_\beta$ where $\beta < \alpha^+$,
and so let
$\beta$ be the least such $\beta.$ Let $Z_\alpha \in \{ 0,1
\}^\alpha \cap L$ be such that $Z_\alpha$ codes $\langle L_\beta,
\in, X_\alpha \rangle$ (we include in $Z_\alpha$ the elementary
diagram of the structure $\langle L_\beta, \in, X_\alpha \rangle)$.
 
For every $\lambda \in E_n \cap B$, let
$$ B_\lambda = \{ \alpha < \lambda: \alpha \in B \; \text{and} \;
Z_\alpha = Z_\lambda | \alpha \}. $$
We have
$$\align B_\lambda &\supseteq \{ \alpha < \lambda : Z_\lambda |
\alpha \; \text{codes} \;
\langle L_\beta, \in, X \rangle \; \text{where} \; \beta \;
\text{is the least} \; \beta \\
&\qquad\text{and} \; X  \; \text{is the least} \; X \; \text{such that}
\; L_\alpha \models \varphi (X) \; \text{and} \; X = X_\lambda \cap
\alpha \} \\
&= \{ \alpha < \lambda : L_\alpha \models \psi (Z_\lambda |
\alpha, X_\lambda \cap \alpha) \} \endalign $$
where $\psi$ is a $\Pi^1_n \wedge \Sigma^1_n$ statement, and
hence
$B_\lambda$ belongs to the filter $F^\lambda_n$.
By the induction hypothesis there is a club
$C_\lambda \subseteq \lambda$ such that $B \cap E_{n-1} \supseteq
B_\lambda \cap E_{n-1} \supseteq C_\lambda \cap E_{n-1}$.
\enddemo
 
\proclaim{Lemma 2.2}  There is a club $C \subseteq \k$ such that $B
\cap E_{n-1}
\supseteq C \cap E_{n-1}$.
\endproclaim
 
\demo{Proof}  If not then $E_{n-1} - B$ is stationary.
This set reflects
at almost all
$\lambda \in E_n$, and since $B \cap E_n$ is stationary, there is
$\lambda \in B \cap
E_n$ such that $(E_{n-1} - B) \cap \lambda$ is stationary in
$\lambda$.  But
$B \cap E_{n-1} \supseteq C_\lambda \cap E_{n-1}$, a
contradiction. \qed\enddemo
 
\definition{Definition 2.3}  For each $t \in L \cap \{ 0,1 \}^{< \k}$, let
$$ S_t = \{ \alpha \in E_{n-1} : t \subset Z_\alpha \}. $$
\enddefinition
 
Since $B \cap E_{n-1}$ is almost all of $E_{n-1}$, there is for
each $\gamma < \k$ some
$t \in \{ 0,1 \}^\gamma$ such that $S_t$ is stationary.
 
\proclaim{Lemma 2.4}  If $t, u \in \{ 0,1 \}^{< \k}$ are such that
both $S_t$ and $S_u$ are stationary
then $t \subseteq u$ or $u \subseteq t$.
\endproclaim
 
\demo{Proof}  Let $\lambda \in B \cap E_n$ be such that
both
$S_t \cap \lambda$ and $S_u \cap \lambda$ are
stationary in $\lambda$.  Let $\alpha, \beta \in C_\lambda$ be
such that
$\alpha \in S_t$ and $\beta \in S_u$.  Since we have $t \subset
Z_\alpha
\subset Z_\lambda$ and
$u \subset Z_\beta \subset Z_\lambda$, it follows that $t
\subseteq u$
or $u \subseteq t$. \qed\enddemo
 
\proclaim{Corollary 2.5}  For each $\gamma < \k$ there is $t_\gamma
\in \{0,1 \}^\gamma$ such that
$S_{t_\gamma}$ is almost all of $E_{n-1}$.
\endproclaim
 
\proclaim{Corollary 2.6}  There is a club $D \subseteq \k$ such that
for all
$\alpha \in D$, if $\alpha \in E_{n-1}$ then
$\alpha \in B$ and $t_\alpha \subset Z_\alpha.$ 
\endproclaim
 
\demo{Proof}  Let $D$ be the intersection of $C$ with the
diagonal intersection of the witnesses for the $S_{t_\gamma}$.\qed
\enddemo
 
\definition{Definition}  $Z = \bigcup \{ t_\gamma : \gamma < \k\}.$
\enddefinition
 
\proclaim{Lemma 2.7}  For almost all $\alpha \in E_{n-1},$ $ Z \cap \alpha
= Z_\alpha$.\endproclaim
 
\demo{Proof}  By Corollary 2.6, if $\alpha \in D \cap E_{n-
1}$ then $Z_\alpha = t_\alpha$.\qed
\enddemo
 
Now we can finish the proof of Lemma 2.1:  The set $Z$  codes a
set $X \subseteq \k$ and witnesses that
$X \in L$.  We claim that $L_\k \models \varphi (X)$.  If not, then
the set $\{ \alpha < \k: L_\alpha \models \neg \varphi (X \cap \alpha )
\}$ is in the filter $F^\k_n$
(because $\neg \varphi$ is $\Sigma^1_n$).
By the induction hypothesis, $L_\alpha \models
\neg \varphi (X \cap \alpha)$
for almost
all $\alpha \in E_{n-1}$.  On the other hand, for almost all
$\alpha \in E_{n-1}$ we
have $L_\alpha \models \varphi (X_\alpha)$ and by Lemma 2.7, for
almost all $\alpha \in E_{n-1}, X \cap \alpha = X_\alpha$; a
contradiction.\qed
\pbf
 
\subhead 3. Proof of Theorem B: Cases 0 and 1
\endsubhead
 
The model is constructed by iterated forcing.  (We refer to [5]
for unexplained notation and terminology).  Iterating with Easton
support, we do a nontrivial construction only at stage $\k$ where
$\k$ is a
inaccessible.
 
Assume that we have constructed the forcing below $\k$, and
denote it $Q$, and denote the model $V (Q)$; if $\lambda < \k$
then $Q | \lambda$ is the forcing below $\lambda$ and
$Q_\lambda \in V (Q | \lambda)$ is the forcing at $\lambda$.  The
rest of
the proof will be to describe $Q_\k$.  The forcing below $\k$ has
size $\k$ and
satisfies the $\k$-chain condition; the forcing at $\k$ will be
essentially $< \k$-closed (for every $\lambda < \k$ has a
$\lambda$-closed dense set) and
will satisfy the $\k^+$-chain condition.  Thus cardinals and
cofinalities are preserved, and stationary subsets of
$\k$ can only be made nonstationary by forcing at $\k$, not below
$\k$ and not after stage $\k$; after stage $\k$ no subsets of
$\k$ are added.
 
By induction, we assume that Full Reflection holds in $V (Q)$ for
subsets of all $\lambda < \k$.  We also assume this for every
$\lambda < \k$:
 
(a)  If $\lambda$ is inaccessible but not weakly compact in $V$
then $\lambda$ is non Mahlo in $V [Q]$.
 
(b)  If $\lambda$ is $\Pi^1_1$-indescribable but not 
$\Pi^1_2$-indescribable in $V$, then
$\lambda$ is 1-Mahlo in $V[Q]$.
 
(c)  And so on accordingly.
 
Let $E_0, E_1, E_2$, etc. denote the subsets of $\k$ consisting
of all inaccessible non Mahlo, 1-Mahlo, 2-Mahlo etc.
cardinals in $V [Q]$.
 
The forcing $Q_\k$ will guarantee Full Reflection for subsets of
$\k$ and make $\k$ into a cardinal of the appropriate Mahlo
class, depending on its indescribability in $V$.  (For instance,
if $\k$ is
$\Pi^1_2$-indescribable but not $\Pi^1_3$-indescribable, it will
be 2-Mahlo in $V (Q \ast Q_\k)$.)
 
The forcing $Q_\k$ is an iteration of length $\k^+$ with 
$< \k$-support of forcing notions that shoot a club through a given set.
We recall ([1], [7], [6]) how one shoots a
club through a single set, and how such forcing iterates:  Given
a
set $B \subseteq \k$, the conditions for shooting a club through
$B$ are closed bounded sets $p$ of ordinals such that $p
\subseteq B$, ordered by end-extension.  In our iteration, the
$B$ will always include the set
Sing of all singular ordinals below $\k$,
which guarantees that the forcing is essentially $< \k$-closed.
One
consequence of this is that at stage $\alpha$ of the iteration,
when
shooting a club through (a name for) a set $B \in V (Q \ast Q_\k
| \alpha)$, the conditions can be taken to be sets in $V (Q)$
rather then
(names for) sets in $V (Q \ast Q_\k | \alpha)$.
 
We use the standard device of iterated forcing:  as $Q_\k$
satisfies the $\k^+$-chain condition, it is possible to enumerate
all names for subsets of $\k$ such that the
$\beta$th name belongs to $V (Q \ast Q_\k | \beta)$, and such
that each name appears cofinally often in the
enumeration.  We call this a {\it canonical enumeration.}
 
We use the following two facts about the forcing:
 
\proclaim{Lemma 3.1}  If we shoot a club through $B$, then every
stationary subset of $B$ remains stationary.
\endproclaim
 
\demo{Proof}  See [5], Lemma 7.38. \qed\enddemo
 
\proclaim{Lemma 3.2}  If $B$ contains a club, then shooting a club
through $B$ has
a dense set that is a $< \k$-closed forcing (and so preserves
all stationary sets).\endproclaim
 
\demo{Proof}  Let $C \subseteq B$ be a club, and let $D =
\{ p: \max (p) \in C \}$. 
\qed\enddemo
 
\remark{Remark}  There is a unique forcing of size $\k$ that is $<
\k$-closed (and nontrivial), namely the one adding a
Cohen subset of $\k$.  We shall henceforth call every forcing
that has such forcing as a dense subset {\it the Cohen forcing} for
$\k$.\endremark
\pmf
 
We shall describe the construction of $Q_\k$ for the cases when
$\k$ is
respectively inaccessible, weakly compact and $\Pi^1_2$-indescribable,
and then outline the general case.  Some details in the three low
cases have to be handled separately from the general case.
\pmf
 
\definition {Case 0}
$Q_\gamma$ for $\gamma$ which is inaccessible
but not weakly compact.\enddefinition
 
We assume that we have constructed $Q | \gamma$, and construct
$Q_\gamma$ in $V (Q | \gamma)$.  To
construct $Q_\gamma$, we first shoot a club through the set $Sing$
and then
do an iteration of length $\gamma^+$ (with $< \gamma$-support),
where at the
stage
$\alpha$ we shoot a club through $B_\alpha$ where $\{ B_\alpha :
\alpha < \gamma^+ \}$ is a
canonical enumeration of all potential subsets of
$\gamma$ such that $B_\alpha \supseteq Sing$.  As $Sing$ contains
a club,
$\gamma$ is in $V (Q \ast P)$ non-Mahlo.  As $Q_\gamma$ is
essentially $< \k$-closed, $\k$ remains inaccessible.
 
In this case, Full Reflection for subsets of $\gamma$ is (vacuously)
true.
\pmf
 
This completes the proof of Case 0.  We shall now introduce some
machinery that (as well as its generalization) we need
later.
 
\definition{Definition 3.3}  Let $\gamma$ be an inaccessible cardinal.
An {\it iteration of order} $0$ (for $\gamma$) is an iteration of
length $< \gamma^+$
such that at each stage $\alpha$ we shoot a club through some
$B_\alpha$
with the property that $B_\alpha \supseteq Sing$.
\enddefinition
 
\proclaim{Lemma 3.4}  (a) If $P$ and $R$ are iterations, and $P$ is of
order $0$ then
$P \force (R$ is of order $0$) if and only of $R$ is of order $0.$
 
(b) If $\dot R$ is a $P$-name then $P \ast \dot R$ is an
iteration of order $0$  if and only if $P$ is an iteration of
order $0$ and
$P \force (\dot R$ is an
iteration of order $0$).
 
(c)  If $A \subseteq Sing$ is stationary and $P$ is an iteration
of order $0$
then $P \force A$ is
stationary.
\endproclaim
\demo{Proof} (a) and (b) are obvious, and (c) is proved as
follows:  Consider the
forcing $R$ that shoots a club through $Sing$.  $R$ is an
iteration
(of length 1) of order $0$, and $R~\ast~P~\force A \text{ is }$
stationary, because $R$ preserves $A$ by Lemma 3.1, and forces
that $P$ is the iterated Cohen forcing
(by Lemma 3.2).  Since $R$ commutes with $P$, we note that $A$ is
stationary in some extension of the forcing extension by $P$, and
so
$P \force A$ is stationary. \qed\enddemo
 
We stated Lemma 3.4 in order to prepare ground for the (less
trivial)
generalization.  We remark that ``$P$ is an iteration of order
$0$'' is a first order property over $V_\gamma$ (using a subset
of
$V_\gamma$ to code the length of the iteration).  The
following lemma, that does not have an analog at higher cases,
simplifies somewhat the
handling of  Case 1.
 
\proclaim{Lemma 3.5}  If $\gamma$ has a $\Pi^1_1$ property $\varphi$ and
$P$ is a $< \gamma$-closed forcing, then $P \force \varphi
(\gamma)$.
\endproclaim
 
\demo{Proof}  Let $\varphi (\gamma) = \forall X \sigma (X)$,
where $\sigma$ is a 1st order
property.  Toward a contradiction, let $p_0 \in P$ and $\dot X$
be such
that $p_0 \force  \neg \sigma (\dot X)$.  Construct a descending
$\gamma$ - sequence of conditions $p_0 \ge p_1 \ge \cdots \ge
p_\alpha \ge
\cdots$ and a continuous sequence $\gamma_0 < \gamma_1 < \cdots <
\gamma_\alpha < \cdots$ such
that for each $\alpha$, $p_\alpha \force \neg \sigma (\dot X
\cap \gamma_\alpha)$, and that $p_\alpha$ decides $\dot X \cap
\gamma_\alpha$; say
$p_\alpha \force \dot X \cap \gamma_\alpha = X_\alpha$.  Let $X
=\bigcup_{\alpha < \gamma} X_\alpha$.
There is a club $C$ such that for all $\alpha \in C,$  $\sigma (X
\cap \alpha)$.  This is a
contradiction since for some $\alpha \in C,$ $\gamma_\alpha =
\alpha$. \qed\enddemo
\pmf
 
\definition{Case 1}  $\lambda$ is $\Pi^1_1$-indescribable but not
$\Pi^1_2$-indescribable.
\enddefinition
 
We assume that
$Q | \lambda$ has been defined, and we shall define an iteration
$Q_\lambda$ of length $\lambda^+$. The idea is to shoot clubs
through the sets
$Sing \cup (Tr (S) \cap E_0)$, for all
stationary sets $S \subseteq Sing$ (including those that appear
at some stage
of the iteration).  Even
though this approach would work in this case, we need to do more
in order to assure that the
construction will work at higher cases.  For that reason we use a
different approach.
 
At each stage of the iteration, we define a filter $F_1$ on
$E_0$, such that the filters all extend the $\Pi^1_1$ filter on
$\lambda$ in $V$, that
the filters get bigger as the iteration progresses, and that sets
that are positive modulo $F_1$ remain positive (and therefore
stationary) at all
later stages.  The iteration consists of shooting clubs through
sets
$B$ such that $B \supseteq Sing$ and $B \cap E_0 \in F_1$, so
that eventually every such $B$ is taken care of.
The crucial property of $F_1$ is that whenever $S$ is a
stationary subset of
$Sing$, then $Tr (S) \cap E_0 \in F_1$.  Thus
at the end of the iteration, every stationary subset of
$Sing$ reflects fully.  Of course, we have to show that the filter
$F_1$ is
nontrivial, that is that in $V (Q | \lambda)$ the set $E_0$ is
positive mod $F_1$.
 
We now give the definition of the filter $F_1$ on $E_0$.  The
definition is nonabsolute enough so that $F_1$ will be
different in each model $V (Q | \lambda \ast Q_\lambda | \alpha)$
for different $\alpha$'s.
 
\definition{Definition 3.6}  Let $C_\lambda$ denote the forcing that
shoots a club through $Sing$.
 
If $\varphi$ is a $\Pi^1_1$ formula and $X \subseteq \lambda$, let
$$ B (\varphi, X) = \{ \gamma \in E_0: \varphi (\gamma, X \cap \gamma )
\} $$
The filter $F_1$ is generated by the sets $B (\varphi, X)$ for those
$\varphi$ and
$X$ such that
$C_\lambda \force \varphi (\lambda, X)$.
A set $A \subset E_0$ is {\it positive} (or {\it 1-positive}), if for
every
$\Pi^1_1$ formula $\varphi$ and every $X \subseteq \lambda$,
if $C_\lambda \force \varphi (\lambda, X)$ then there exists a
$\gamma \in A$ such that $\varphi (\gamma, X \cap \gamma)$.
\enddefinition
 
\remark{Remarks}\endremark
 
1.  The filter $F_1$ extends the club filter (which is generated
by the sets
$B (\varphi, X)$ where $\varphi$ is first-order).  Hence every positive
set is stationary.
 
2.  The property ``$A$ is 1-positive'' is $\Pi^1_2$.
 
\proclaim{Lemma 3.7}  In $V (Q | \lambda),$ $E_0$ is positive.
\endproclaim
 
\demo{Proof}  We recall that in $V,$ $\lambda$ is $\Pi^1_1$-indescribable, 
and $E_0$ is the set of inaccessible, non-weakly-compact
cardinals.  Let $Q = Q | \lambda$.  So let $\varphi$ be a $\Pi^1_1$
formula,
let $\dot X$ be a $Q$ -
name for a subset of $\lambda$, and assume that $V (Q \ast
C_\lambda) \models
\varphi (\lambda, \dot X)$. The statement that $Q \ast C_\lambda
\force \varphi
(\lambda, \dot X)$ is a $\Pi^1_1$ statement (about $Q,$ $C$ and
$\dot X$).  By
$\Pi^1_1$-indescribability, this reflects to some $\gamma \in
E_0$ (as $E_0$ is positive in the
$\Pi^1_1$ filter).  Since $Q \cap V_\gamma = Q | \gamma$ and
since $Q | \gamma$ satisfies the $\gamma$-chain condition, the
name $\dot X$ reflects to the $Q | \gamma$-name for $\dot X
\cap \gamma$.
Also $C_\lambda \cap V_\gamma = C_\gamma$.  Hence
$$ Q | \gamma \ast C_\gamma \force \varphi (\gamma, \dot X \cap
\gamma). $$
What we want to show is that $V (Q) \models \varphi (\gamma, \dot X
\cap \gamma)$.
Since forcing above $\gamma$ does not add subsets of $\gamma$ it
is enough to show that
$V (Q | \gamma \ast Q_\gamma) \models \varphi$.  However, $C_\gamma$
was the first stage of $Q_\gamma$ (see Case 0), and the rest of
$Q_\gamma$ is the
iterated Cohen forcing for $\gamma$.  By Lemma 3.5, if $\varphi$ is
true in $V (Q | \gamma \ast C_\gamma)$, then it is true in $V (Q
| \gamma \ast Q_\gamma)$. \qed\enddemo
 
\proclaim{Lemma 3.8} If $S \subseteq Sing$ is stationary, then the set
$\{ \gamma \in E_0: S \cap \gamma$ is stationary$\}$ is in
$F_1$.
\endproclaim
 
\demo{Proof}  The property $\varphi (\lambda, S)$ which states
that $S$ is stationary is $\Pi^1_1$.  If we show that $C_\lambda
\force \varphi (\lambda, S)$,
then $\{ \gamma \in E_0: \varphi (\gamma, S \cap \gamma ) \}$ is in
$F_1$.
But forcing with $C_\lambda$ preserves stationarity of $S$, by
Lemma 3.1. \qed\enddemo
 
\definition{Definition 3.9}  An {\it iteration of order 1} (for
$\lambda$) is an
iteration of length $< \lambda^+$ such that at each stage
$\alpha$ we shoot
a club through some $B_\alpha$ such that $B_\alpha \supseteq
Sing$ and
$B_\alpha \cap E_0 \in F_1$.
\enddefinition
 
\remark{Remark} If we include the witnesses for $B_\alpha \cap E_0
\in F_1$ as parameters in
the definition, i.e. $\varphi_\alpha,$ $X_\alpha$ such that $C_\lambda
\force \varphi_\alpha (\lambda, X_\alpha)$ and $B_\alpha
\cap E_0 \supseteq \{ \gamma \in E_0: \varphi (\gamma, X_\alpha \cap
\gamma ) \}$, then the property ``$P$ is an iteration of order 1"
is $\Pi^1_1$.  \endremark
\pmf
 
We shall
now give the definition of $Q_\lambda$:
\definition{Definition 3.10}  $Q_\lambda$ is (in $V (Q (\lambda))$ an
iteration of length $\lambda^+$, such that for each $\alpha <
\lambda^+,$ $Q_\lambda | \alpha$ is an iteration of
order 1, and such that each potential $B$ is used as $B_\beta$ at
cofinally many stages  $\beta$.
\enddefinition
 
We will now show that both ``$B \in F_1$" and ``$A$ is positive"
are preserved under
iterations of order 1:
 
\proclaim{Lemma 3.11}  If $B \in F_1$ and $P$  is an iteration of
order 1 then $P \force B \in F_1$.  Moreover, if
$(\varphi, X)$ is a witness for $B \in F_1$, then it remains a
witness after forcing with $P$.
\endproclaim
 
\demo{Proof}  Let $B \supseteq B (\varphi, X)$ where $\varphi$ is
$\Pi^1_1$
and $C_\lambda \force \varphi (\lambda, X)$, and let $P$ be an
iteration of order 1.  As $P$ does not
add bounded subsets, $B (\varphi, X)$ remains the same, and so we
have to
verify that $P \force (C_\lambda \force \varphi)$.  However,
$C_\lambda$ commutes with $P$,
and moreover, $C_\lambda$ forces that $P$ is the Cohen forcing
(because
after $C_\lambda,$ $P$ shoots clubs through sets that contain a
club,
see Lemma 3.2). By Lemma 3.5, $C_\lambda
\force \varphi$ implies that $C_\lambda \force ( P \force \varphi)$.
\qed\enddemo
 
\proclaim{Lemma 3.12}  If $A \subseteq E_0$ is positive and $P$ is an
iteration of order 1 then $P \force A$ is positive.
\endproclaim
 
We postpone the proof of this crucial lemma for a while.  We
remark that the
assumption under which Lemma 3.12 will be proved is that the
model in which we are working contains $V (Q | \lambda)$; this
assumption will be satisfied in the future when the Lemma is
applied.
 
\proclaim{Lemma 3.13} (a)  If $P$ and $R$ are iterations, and $P$ is of
order 1 then
$P \force (R$ is of order 1) if and only if $R$ is of order 1.
 
(b) If $\dot R$ is a $P$-name  then $P \ast \dot R$ is an
iteration of order 1 if and only if $P$ is an iteration of order
1 and $P \force (\dot R$ is an iteration of order 1).
 
(c) Every iteration of order 1 is an iteration of order $0$.
\endproclaim
 
\demo{Proof}  Both (a) and (b) are consequences of Lemma
3.12.  The decision whether a particular stage of the iteration
$R$
satisfies the definition of being of order 1 depends only on whether
$B_\alpha \in F_1$, which does not depend on $P$.
 
(c) is trivial. \qed\enddemo
 
\proclaim{Corollary 3.14}  In $V (Q | \lambda \ast Q_\lambda), E_0$ is
stationary (so
$\lambda$ is 1-Mahlo), and every stationary $S \subseteq Sing$
reflects fully in
$E_0$.
\endproclaim
 
\demo{Proof}  Suppose that $E_0$ is not stationary.  Then
it is disjoint
from some club $C$, which appears at some stage $\alpha <
\lambda^+$ of the
iteration $Q_\lambda$.  So $E_0$ is nonstationary in $V (Q |
\lambda \ast
Q_\lambda | (\alpha +1))$.  This is a contradiction, since $E_0$
is positive in that model, by Lemmas 3.7 and 3.12.
 
If $S$ is a stationary subset of
$Sing$, then $S \in V (Q | \lambda \ast Q_\lambda | \alpha)$ for
some
$\alpha$ and so by Lemma 3.8, $B = Tr (S) \cap E_0 \in F_1$ (in
that model).
Hence $B$ remains in $F_1$ at all later stages, and eventually,
$B = B_\alpha$ is used at stage $\alpha$, that is we produce a
club $C$ so that $B \supseteq
C \cap E_0$.  Since $Q_\lambda$ adds no bounded subsets of
$\lambda$, the trace of $S$ remains the same, and so $S$ reflects
fully in $V (Q | \lambda \ast Q_\lambda)$. \qed\enddemo
 
\subhead{Proof of Lemma 3.12}\endsubhead
 
Let $\varphi$ be a $\Pi^1_1$ property, and let $\dot X$ be a $P$-name 
for a subset of
$\lambda$.  Let $p \in P$ be a condition that forces that
$C_\lambda \force \varphi (\lambda, \dot X)$.
We are going to find a stronger $q \in P$ and a $\gamma \in A$
such that $q$ forces $\varphi (\gamma, \dot X \cap \gamma).$
 
$P$ is an iteration of order 1, of length $\alpha$.  At stage
$\beta$ of the iteration, we have
$P | \beta$ - names $\dot B_\beta,$ $\varphi_\beta$ and $\dot X_\beta$
for a set
$\supseteq Sing$, a $\Pi^1_1$ formula, and a subset of $\lambda$
such that $P | \beta$ forces
that $C_\lambda \force \varphi_\beta (\lambda, \dot X_\beta)$ and
that $\dot B_\beta \supseteq \{ \gamma \in E_0: \varphi_\beta
(\gamma, \dot X_\beta \cap \gamma ) \}$, and we shoot a club
through $\dot B_\beta$.
 
Let $\psi$ be the following statement (about $V_\lambda$ and a
relation on $V_\lambda$ that codes a model of size $\lambda$
including the relevant parameters and satisfying enough
axioms of ZFC; the relation will also insure that the model $M$
below has the properties that we list):
\pmf
 
{\sl $P$ is an iteration of length $\alpha$, at each stage shooting a
club through $\dot B_\beta \supseteq Sing$,
and $p \force \varphi (\lambda, \dot X)$ and for every
$\beta < \alpha,$ $P | \beta \force \varphi_\beta (\lambda, \dot
X_\beta)$.}
\pmf
 
First we note that $\psi$ is a $\Pi^1_1$ property.  Secondly, we
claim that
$C_\lambda \force \psi$:  In the forcing extension by $C_\lambda,$
$P$ is still an iteration etc., and $p \force \varphi$ and $P | \beta
\force \varphi_\beta$
because in the ground model, $p \force (C_\lambda \force \varphi)$
and
$P | \beta \force (C_\lambda \force \varphi_\beta)$,
and $C_\lambda$ commutes with $P$.
 
Thus, since $A$ is positive in the ground model, there exists
some $\gamma \in A$ such that $\psi (\gamma,\text{ parameters}\cap
V_\gamma$).  This
gives us a model $M$ of size $\gamma$, and its transitive
collapse $N = \pi (M)$, with the following properties:
 
(a)  $M \cap \lambda = \gamma$,
 
(b) $P, p, \dot X \in M$ and $M \models P$ is an iteration given
by $\{ \dot B_\beta : \beta < \alpha \}$,
 
(c) $p \force \varphi (\gamma, \pi (\dot X))\qquad$  (the forcing $\force$
is in $\pi (P)$),
 
(d)  $\forall \beta < \alpha$, if $\beta \in M$, then $\pi (P |
\beta) \force \varphi_\beta (\gamma, \pi (\dot X_\beta))$.
\pmf
It follows that $\pi (P)$ is an iteration on $\gamma$ (or order
$0$), of length $\pi (\alpha)$,
that at stage $\pi (\beta)$ shoots a club through
$\pi (\dot B_\beta)$.  Also, $p \force \pi (\dot X) = \dot X \cap
\gamma$ (forcing in $P$).
 
\proclaim{Sublemma 3.12.1}  There exists an $N$-generic filter $G \ni
p$ on $\pi (P)$ such that if $X \subseteq \gamma$ denotes the $G$-interpretation
$\pi (\dot X) / G$ of $\pi (\dot X)$, and for each $\beta \in M,
X_\beta = \pi (\dot X_\beta ) / G$, then $\varphi (\gamma, X)$ and
$\varphi_\beta (\gamma, X_\beta)$ hold.
\endproclaim
 
\demo{Proof}  We assume that $V (Q | \lambda)$ is a part of
our
universe, and that no
subsets of $\gamma$ have been added after $Q_\gamma$.  So it
suffices to find $G$ in
$V (Q | \gamma \ast Q_\gamma)$.  Note also
that $E_0 \cap \gamma$ is nonstationary (as $\gamma$ was made non
Mahlo by $Q_\gamma$).
Since $\pi (P)$ is an
iteration of order $0$, since $Sing$ contains a club, and because
$\pi (P)$ has size $\gamma$,
it is the Cohen forcing for $\gamma$,
and therefore isomorphic to the forcing at each stage of the
iteration
$Q_\gamma$ except the first one (which is $C_\gamma$).
 
There is $\eta < \gamma^+$ such that $V (Q | \gamma \ast Q_\gamma
| \eta)$ contains
$\pi (P), \pi (\dot X)$, all members
of $N$, and all $\pi (\dot X_\beta), \beta \in M$.  Also, the
statements $p \force
\varphi (\gamma, \pi (\dot X))$ and $\pi (P | \beta)
\force \varphi_\beta (\gamma, \pi (\dot X_\beta))$,
being $\Pi^1_1$ and true, are true in $V (Q | \gamma \ast
Q_\gamma | \eta)$.
As $\pi (P)$ (below $p$) as
well as the $\pi (P | \beta)$ are isomorphic to the $\eta$$^{th}$
stage $Q_\gamma (\eta)$ of $Q_\gamma$, and we do have a generic
filter
for $Q_\gamma (\eta)$ over $V (Q | \gamma \ast Q_\gamma | \eta)$,
we have a $G$ that is
$N$-generic for $\pi (P)$ and $\pi (P | \beta)$.  If we let
$X = \pi (\dot X) /
G$ and $X_\beta = \pi (\dot X) / G$, then in $V (Q | \gamma \ast
Q_\gamma | (\eta +1))$ we have $\varphi (\gamma, X)$ and $\varphi_\beta
(\gamma, X_\beta)$.
Since the rest of the iteration $Q_\gamma$ is the iterated Cohen
forcing, we
use Lemma 3.5 again to conclude that $\varphi (\gamma, X)$ and
$\varphi_\beta (\gamma, X_\beta)$ are true
in $V (Q | \gamma \ast Q_\gamma)$, hence are true. \qed\enddemo
 
Now let $H = \pi^{-1} (G)$ and for every $\beta \in M$ let
$B_\beta = \pi (\dot B_\beta) / G$.
By induction on $\beta \in M$, we construct a condition $q \le p$
(with support
$\subseteq M$) as follows:  For each $\xi \in M$, let $q (\xi) =
H_\xi \cup \{ \lambda \}$.  This is a closed set of ordinals.  At
stage $\beta,$ $q | \beta$ a condition by the induction hypothesis, and
$q | \beta \supseteq H | \beta$ (consequently, $q | \beta$ forces
$\dot X_\beta \cap \gamma = X_\beta$ and $\dot B_\beta \cap
\gamma = B_\beta$).
$H_\beta$ is a closed set of ordinals,
cofinal in $\gamma$, and $H_\beta \subseteq B_\beta$.  We let $q
(\beta) = H_\beta \cup \{ \gamma \}$.  In order that
$q | (\beta +1)$ is a condition it is necessary that $q | \beta
\force
\gamma \in \dot B_\beta$.  But by Sublemma 3.12.1 we have
$\varphi_\beta
(\gamma, X_\beta)$, so this is forced by $P$ (which
does not add subsets of $\gamma$), and since $q | \beta \force
X_\beta =
\dot X_\beta \cap \gamma$, we have
$q | \beta \force \varphi_\beta (\gamma, \dot X_\beta \cap \gamma)$.
But this implies that
$q | \beta \force \gamma \in \dot B_\beta$.  Hence $q | (\beta
+1)$ is a condition, which extends
$H | (\beta +1)$.
 
Therefore $q$ is a condition, and since $q \supseteq H$, we have
$q \force
\dot X \cap \gamma = X$.
But $\varphi (\gamma, X)$ holds by Sublemma 3.12.1., so it is forced
by $q$, and
so $q \force \varphi (\gamma, \dot X \cap \gamma)$, as
required. 
\qed
\pbf
 
\subhead 4. Case 2 and up \endsubhead
 
Let $\k$ be $\Pi^1_2$-indescribable but not $\Pi^1_3$-indescribable.  
Below $\k$,
we have four different types of limit cardinals in $V$:
 
$Sing =$ the singular cardinals
 
$E_0 =$ inaccessible not weakly compact
 
$E_1 = \Pi^1_1$- but not $\Pi^1_2$-indescribable
 
the rest $= \Pi^1_2$ indescribable
\pmf
We shall prove a sequence of lemmas (and give a sequence of
definitions), analogous to
3.6--3.14.  Whenever possible, we use the same argument; however,
there are some changes and additional
complications.
 
\definition{Definition 4.1}  A $\Pi^1_2$ formula $\varphi$ is {\it absolute}
for $\lambda \in E_1$ if for every
$\alpha < \lambda^+$ and every $X \in V (Q | \lambda \ast
Q_\lambda | \alpha)$,
\pmf
(1) $V (Q | \lambda \ast Q_\lambda | \alpha) \models$ (for every
iteration $R$ of order 1, $\varphi (\lambda, X)$ iff
$R \force \varphi (\lambda, X)$),
\pmf
(2) $V (Q | \lambda \ast Q_\lambda | \alpha) \models \varphi
(\lambda, X)$ implies $V (Q | \lambda \ast Q_\lambda) \models
\varphi
(\lambda, X)$, and
\pmf
(3) $V (Q | \lambda \ast Q_\lambda | \alpha) \models \neg \varphi
(\lambda, X)$ implies $V (Q | \lambda \ast Q_\lambda) \models
\neg \varphi
(\lambda, X)$.
\par\bigpagebreak
We say that $\varphi$ is {\it absolute} if it is absolute for all
$\lambda \in E_1,$ $\lambda < \k$.
\enddefinition
 
\definition{Definition 4.2}  If $\varphi$ is a $\Pi^1_2$ formula and $X
\subseteq \k$, let
$$ B (\varphi, X) = \{ \lambda \in E_1  : \varphi (\lambda, X \cap
\lambda) \}.$$
The filter $F_2$ is generated by the sets $B (\varphi, X)$ where $\varphi$
is an absolute $\Pi^1_2$
formula and $X$ is such that
$R \force \varphi (\k, X),$ for all iterations $R$ of order 1.
 
A set $A \subseteq E_1$ is {\it positive}
(2-positive) if for any absolute $\Pi^1_2$ formula $\varphi$ and
every $X \subseteq \k$, if every
iteration $R$ of order 1 forces $\varphi (\k, X)$, then there exists
a $\lambda \in A$ such that $\varphi (\lambda, X \cap \lambda)$.
\enddefinition
 
\remark{Remark}  The property ``A is 2-positive" is $\Pi^1_3$.
\endremark
 
\proclaim{Lemma 4.3}  In $V (Q | \k),$ $E_1$ is positive.
\endproclaim
 
\demo{Proof}  Let $Q = Q | \k$.  Let $\varphi$ be an absolute
$\Pi^1_2$ formula, and
let $\dot X$ be a $Q$-name for a subset of $\k$, and assume that
in $V (Q),$ $R \force \varphi (\k, \dot X)$ for all order-1
iterations $R$.  In particular, (taking $R$ the empty iteration),
$V (Q) \models \varphi (\k, \dot X)$.
 
Using the $\Pi^1_2$-indescribability of $\k$ in $V$, there exists
a $\lambda \in E_1$
such that $V (Q | \lambda) \models \varphi (\lambda, \dot X \cap
\lambda)$.  In order
to prove that $V (Q) \models \varphi (\lambda, \dot X \cap
\lambda)$, it is enough to show that
$V (Q | \lambda \ast Q_\lambda) \models \varphi (\lambda, \dot X
\cap \lambda)$.
This however is true because $\varphi$ is absolute for $\lambda$.\qed
\enddemo
 
\proclaim{Lemma 4.4}  The property ``S is 1-positive" of a set $S
\subseteq E_0$ is an absolute $\Pi^1_2$ property,
and is preserved under forcing with iterations of order 1.
\endproclaim
 
\demo{Proof} The preservation of ``1-positive" under
iterations of order 1 was proved in Lemma 3.12.
To show that the property is absolute for all $\lambda \in E_1$,
first assume
that $S \in V (Q | \lambda \ast Q_\lambda | \alpha)$ is 1-positive.  Since all
longer initial segments of the iteration $Q_\lambda$ are
iterations of order 1, hence order 1 iterations over
$Q_\lambda | \alpha$ (by Lemma 3.13), $S$ is 1-positive in each
$V (Q | \lambda \ast Q_\lambda | \beta),$ $\beta > \alpha$.
However, the property ``S is 1-positive" is $\Pi^1_2$, and so it
also holds in
$V (Q | \lambda \ast Q_\lambda)$, because every subset of
$\lambda$ in
that model appears at some stage $\beta$.  (We remark that this
argument, using $\Pi^1_2$,
does not work in higher cases).
 
Conversely, assume that $S$ is not 1-positive in $V (Q |
\lambda \ast Q_\lambda | \alpha)$.  There exists a $\Pi^1_1$
formula $\varphi$ and some
$X \subseteq \lambda$ such that $\varphi (\gamma, X \cap \gamma)$
fails for all $\gamma \in S$,
while $C_\lambda \force \varphi (\lambda, X)$.  The rest of the
argument is the same as the one in Lemma 3.11:
Let $P = Q_\lambda / (Q_\lambda | \alpha);$ $C_\lambda$ commutes
with $P$ and forces that $P$
is the iterated Cohen forcing.  Hence by Lemma 3.5, $P \force
(C_\lambda \force \varphi)$, i.e. $V (Q | \lambda \ast
Q_\lambda) \models
(C_\lambda \force \varphi)$.  Therefore $S$ is not 1-positive in $V
(Q | \lambda \ast Q_\lambda)$.
(Again, this argument does not work in higher cases.). \qed\enddemo
 
\proclaim{Lemma 4.5}  The property ``R is an iteration of order 1" is
an absolute $\Pi^1_2$ property,
and is preserved under forcing with iterations of order 1.
Moreover, in $V (Q | \lambda \ast Q_\lambda)$, if $R$ is an
iteration
of order 1, then $R$ is the Cohen forcing.
\endproclaim
 
\demo{Proof}  The preservation of the property under
iterations of order 1 was
proved in Lemma 3.13.  If $R$ is an iteration of order 1 in $V (Q
| \lambda \ast
Q_\lambda | \alpha)$, shooting clubs through $\dot A_0, \dot A_1,
\dot A_2$, etc., then $R$
embeds in $Q_\lambda$ above $\alpha$ as a subiteration, i.e.
there are $\beta_0, \beta_1$, etc. such that
$\dot A_0 = B_{\beta_0}, \dot A_1 = B_{\beta_1}$, etc.  Moreover,
there is some $\gamma > \alpha$ such that the
$A_0, A_1, A_2$, etc. all contain a club.  Hence $R$ is the Cohen
forcing in $V (Q | \lambda \ast Q_\lambda | \gamma)$.
Therefore $R$ is the Cohen forcing in $V (Q | \lambda \ast
Q_\lambda)$, and
consequently an iteration of order 1.  As for the absoluteness
downward, we give the proof for
iterations of length 2.  Let $M_\infty = V (Q | \lambda \ast
Q_\lambda)$, let $R = (R_0, R_1)$ be an iteration given by
$A_0$ and $\dot A_1 \in M_\infty (R_0)$, such that 
in $M_\infty,$ $A_0 \in F_1$ and
$R_0 \force \dot A_1 \in F_1$.  Let $R \in M_\alpha = V (Q |
\lambda \ast Q_\lambda | \alpha)$.
We will show that in $M_\alpha,$ $R$ is an iteration of order 1,
and that in $M_\infty,$ $R$ is the Cohen forcing.
 
First, since $A_0 \in F_1$ is absolute, there is a $\beta >
\alpha$ such that
$M_\beta \models A_0$ contains a club and such that $A_0 =
B_\beta$ ($B_\beta$ is the set
used at stage $\beta$ of the iteration $Q_\lambda$).  Since
$M_\beta \models
(R_0$ is Cohen), we have $M_\infty \models R_0$ is Cohen.
 
Now, in $M_\infty$ we have $R_0 \force \dot A_1 \in F_1$.  We
claim that
in $M_\beta,$ $R_0 \force \dot A_1 \in F_1$.   Then it follows that
$R$ is an iteration of order 1 in $M_\beta$.
 
It remains to prove the claim.  Let $\dot X$ denote $\dot A_1$,
let $\varphi (\dot X)$ denote the
absolute $\Pi^1_2$ property $\dot A_1 \in F_1$ and let $C$ denote
the Cohen forcing.  We recall
that $M_{\beta +1} = M_\beta (C)$.\enddemo
 
\proclaim{Sublemma 4.5.1}  Let $\dot X$ be a $C$-name in $M_\beta$,
and assume that
$M_{\beta +1} = M_\beta (C)$.  If $C \force \varphi
(\dot X)$ in $M_\infty,$ then
$C \force \varphi (\dot X)$ in $M_\beta$.
\endproclaim
 
\demo{Proof}  Let $P$ be the forcing such that $M_\infty =
M_{\beta +1} (P)$, and assume, toward a contradiction,
that $C \force \varphi (\dot X)$ in $M_\infty$ but $C \force \neg
\varphi (\dot X)$ in $M_\beta$.  Let
$G_C \times G_P \times H$ be a generic on $C \ast P \ast C$, and
let $X = \dot X / H$.  Let $C = C_1 \times C_2$
where both $C_1$ and $C_2$ are Cohen, and consider the generic $H
\times G_C \times G_P$ on
$C_1 \times C_2 \times P = C \times P$ (it is a generic
because since $H$ is generic over $G_C \times G_P,$ $G_C \times
G_P$ is generic over $H$).
 
In $M_\beta,$ $C_1$ forces $\varphi$ false, hence $\varphi (X)$ is false
in $M_\beta [H]$.
Since $\neg \varphi$ is preserved by Cohen forcing (in fact by all
order-1 iterations), so $\varphi (X)$ is false in $M_\beta [H \times
G_C]$.
Now $\varphi$ is absolute (between $M_{\beta +1}$ and $M_\infty$)
and so $\varphi (X)$ is false in
$M_\beta [H \times G_C \times G_P]$.  On the other hand, since $C
\force \varphi (\dot X)$ in $M$, we have
$M_\beta [G_C \times G_P \times H] \models \varphi (\dot X / H)$, so
$\varphi$ is
true in $M_\beta [G_C \times G_P \times H]$, a contradiction.
\qed\enddemo
 
\proclaim{Lemma 4.6}  If $S \subseteq E_0$ is 1-positive, then the set
$$ \{ \lambda \in E_1 : S \cap \lambda \; \text{is 1-positive} \}
$$
is in $F_2$.  Therefore $Tr (S) \cap E_1 \in F_2$.
\endproclaim
 
\demo{Proof}  The first sentence follows from the
definition of $F_2$ because ``$S$ is 1-positive"
is absolute $\Pi^1_2$ and if $S$ is positive then it is positive
after every order 1 iteration.  The
second sentence follows, since 1-positive subsets of $\lambda$
are stationary. 
\qed\enddemo
 
\definition{Definition 4.7}  An {\it iteration of order 2} (for $\k$) is
an iteration of length $< \k^+$ that at
each stage $\alpha$ shoots a club through some $B_\alpha$ such
that $B_\alpha
\supseteq Sing,$ $B_\alpha \cap E_0 \in F_1$, and $B_\alpha \cap
E_1 \in F_2$.
\enddefinition
 
\remark{Remarks}  1.  An iteration of order 2 is an iteration of
order 1.
 
2. If we include the witnesses for $B_\alpha$ to be in the
filters, then the property ``$P$ is an iteration of order 2" is
$\Pi^1_3$.
\endremark
 
\definition{Definition 4.8}  $Q_\k$ is (in $V (Q | \k)$) an iteration of
length $\k^+$, such
that for each $\alpha < \k^+,$ $Q_\k | \alpha$ is an iteration of
order 2, and such that each potential $B$ is
used as $B_\beta$ at cofinally many stages $\beta$.
\enddefinition
 
\proclaim{Lemma 4.9}  If $B \in F_2$ and $P$ is an iteration of order
2 then $P \force B \in F_2$ (and a witness $(\varphi, X)$
remains a witness).
\endproclaim
 
\demo{Proof}  Let $B \supseteq B (\varphi, X)$ where $\varphi$ is
an absolute $\Pi^1_2$,
and every iteration of order 1 forces $\varphi$; let $P$ be an
iteration of order 2.  Since $P$
does not add subset of $\k,$ $B (\varphi, X)$ remains the same and
$\varphi$ remains absolute.
Thus it suffices to verify that for each $P$-name $\dot R$ for
an order 1
iteration, $P \force (\dot R \force \varphi)$.  However, $P$ is an
iteration of order 1, so by
Lemma 3.13, $P \ast \dot R$ is an iteration of order 1, and by
the assumption on $\varphi,$ $P \ast \dot R \force \varphi$. 
\qed\enddemo
 
\proclaim{Lemma 4.10}  If $A \subseteq E_1$ is 2-positive and $P$ is
an iteration of order 2 then $P \force A$ is
2-positive.
\endproclaim
 
\demo{Proof}  Let $\varphi$ be an absolute $\Pi^1_2$ property,
let $\dot X$ be a $P$-name for a subset of $\k$ and let
$p \in P$ force that for all order-1-iterations $R,$ $R \force \varphi
(\k, \dot X)$.  We want a $q \le p$ and a $\lambda \in A$ such
that $q \force \varphi (\lambda, \dot X \cap \lambda)$.
 
$P$ is an iteration of order 2 that at each stage $\beta$ (less
than the length
of $P$) shoots a club through a set $\dot B_\beta$ such that $P |
\beta$ forces
that
 
(1) $\dot B_\beta \supseteq Sing$,
 
(2) $\dot B_\beta \cap E_0 \supseteq \{ \gamma \in E_0:
\varphi^1_\beta (\gamma, \dot X_\beta \cap \gamma)\}$, and
 
(3) $\dot B_\beta \cap E_1 \supseteq \{ \lambda \in E_1:
\varphi^2_\beta (\lambda, \dot Y_\beta \cap \lambda)\}$,
 
where $\dot X_\beta$ and $\dot Y_\beta$ are names for subsets of
$\k$, the $\varphi^1_\beta$ are
$\Pi^1_1$ formulas (with some extra property that make $P$ an
order-1 interation) and
the $\varphi^2_\beta$ are absolute (in $V (Q | \k))$ $\Pi^1_2$
properties,
and $P | \beta$ forces that $\forall R$ (if $R$ is an iteration
of order
1 then $R \force \varphi^2_\beta (\k, \dot Y_\beta))$.
 
We shall reflect, to some $\lambda \in E_1$, the $\Pi^1_2$
statement $\psi$ that states (in addition to a
first order statement in some parameter that produces the model
$M$ below):
 
(a) $P$ is an iteration of order 1 using the $\varphi^1_\beta, \dot
X_\beta, \varphi^2_\beta, \dot Y_\beta$,
 
(b) $p \force \varphi (\k, \dot X)$,
 
(c) for every $\beta <$ length$(P),$ $P | \beta \force
\varphi^2_\beta (\k, \dot Y_\beta)$.
 
First we note that $\psi$ is a $\Pi^1_2$ property.  Secondly, we
claim that
$\psi$ is absolute for every $\lambda \in E_1$.  Being an
iteration of order 1 is
absolute by Lemma 4.5.  That (b) and (c) are absolute will follow
once we show that if $\varphi$ is an absolute $\Pi^1_2$
property and $R$ an iteration of order 1, then ``$R \force \varphi$"
is absolute:
\enddemo
 
\proclaim{Sublemma 4.10.1}.  Let $\varphi$ be absolute for $\lambda$, let
$\alpha < \lambda^+,$ $X, R \in V (Q | \lambda \ast Q_\lambda |
\alpha)$ be a subset of
$\lambda$ and an iteration of order 1.  Then the property $R
\force \varphi (\lambda, X)$ is absolute between $M_\alpha = V (Q |
\lambda \ast Q_\lambda | \alpha)$ and $M_\infty =
V (Q | \lambda \ast Q_\lambda)$.
\endproclaim
 
\demo{Proof}  Let $M_\alpha \models (R \force \varphi
(\lambda, X))$.  Then
$M_\alpha \models \varphi (\lambda, X)$ and by absoluteness,
$M_\infty \models \varphi (\lambda, X)$.  If
in $M_\infty,$ $R \force \neg \varphi (\lambda, X)$, then because $R$
is in $M_\infty$ the Cohen forcing, there is (by Sublemma 4.5.1)
some $\beta > \alpha$ such
that $R$ is the Cohen forcing in $M_\beta$ and $M_\beta \models
(R \force \neg \varphi)$.  By absoluteness again, $M_\beta \models
\neg \varphi$, a contradiction.
\qed\enddemo
 
Thus $\psi$ is an absolute $\Pi^1_2$ property.  Next we show that
if $R$ is an iteration of order 1
then $R$ forces $\psi (\k$, parameters):
 
(a)  $R \force (P$ is an iteration of order 1), by Lemma 3.13.
 
(b) $R$ commutes with $P$, and by the assumption of the proof, $p
\force (R \force \varphi (\k, \dot X))$.  Hence
$R \force (p \force \varphi (\k, \dot X))$.
 
(c)  For every $\beta,$ $R$ commutes with $P | \beta$, and by the
assumption on $\varphi^2_\beta,
$ \newline $
P | \beta \force (R \force
\varphi^2_\beta (\k, \dot Y_\beta))$. Hence
$R \force (P | \beta \force \varphi^2_\beta (\k, \dot Y_\beta))$.
 
Now since $A$ is 2-positive in the ground model, there exists a
$\lambda \in A$ such
that \newline $\psi (\lambda,\text{ parameters} \cap V_\lambda$).  This gives
us a model $M$ of
size $\lambda$, and its transitive collapse $N = \pi (M)$, with
the following properties:
 
(a)  $M \cap \k = \lambda$,
 
(b)  $P, p, \dot X \in M$,
 
(c) $\pi (P)$ is an iteration of order 1 for $\lambda$,
 
(d)  $p \force \varphi (\lambda, \pi (\dot X))$,
 
(e) $\forall \beta \in M \quad \pi (P | \beta) \force
\varphi^2_\beta (\lambda,
\pi (\dot Y_\beta))$.
 
The rest of the proof is analogous to the proof of Lemma 3.12, as
long as
we prove the analog of Sublemma 3.12.1:  after that, the proof is
Case 1 generalizes
with the obvious changes.
 
\proclaim{Sublemma 4.10.2} There exists an $N$-generic filter $G \ni
p$ on $\pi (P)$ such that if $X = \pi (\dot X) / G$
and $Y_\beta = \pi (\dot Y_\beta)/ G$ for each $\beta \in M$,
then $\varphi (\lambda, X)$ and
$\varphi^2_\beta (\lambda, Y_\beta)$ hold.
\endproclaim
 
\demo{Proof}  We find $G$ in $V (Q | \lambda \ast
Q_\lambda)$.  Since $\pi (P)$ is
an iteration of order 1 and
$\varphi$ is absolute, there is an $\alpha < \lambda^+$ such that $V
(Q | \lambda \ast Q_\lambda | \alpha)$ contains $\pi (\dot X),
\pi (\dot Y_\beta)$ ($\beta \in M)$ and the dense sets in $N$,
thinks that $\pi (P)$ is the Cohen
forcing, such that the forcing $Q_\lambda (\alpha)$ at stage
$\alpha$ is the Cohen forcing, and (by absoluteness  and by
Sublemma 4.10.1) $Q_\lambda (\alpha)$ (or $\pi (P)$) forces
$\varphi (\lambda, \pi (\dot X))$ and $\varphi^2_\beta (\lambda, \pi
(\dot Y_\beta))$.  The generic filter on
$Q_\lambda (\alpha)$ yields a generic $G$ such that $V (Q |
\lambda \ast Q_\lambda | (\alpha +1)) \models \varphi (\lambda, X)$
and $\varphi^2_\beta (\lambda, Y_\beta)$ where $X = \pi (\dot X) /
G,$ $Y_\beta = \pi (\dot Y_\beta) / G$.  By absoluteness again,
$\varphi (\lambda, X)$ and $\varphi^2_\beta(\lambda, Y_\beta)$
hold in $V (Q | \lambda \ast Q_\lambda)$, and hence they hold.
\qed\enddemo
 
\proclaim{Lemma 4.11}  (a)  If $P$ and $R$ are iterations, and $P$ is
of order 2, then $P \force (R$ is of order 2) if and only if $R$
is order 2.
 
(b)  If $\dot R$ is a $P$-name then $P \ast \dot R$ is an
iteration of order 2 if and only if $P$ is an iteration of order
2 and
$P \force (\dot R$ is an iteration of order 2).
\endproclaim
 
\demo{Proof}  By Lemma 4.10 (just as Lemma 3.13 follows
form Lemma 3.12). 
\qed\enddemo
 
\proclaim{Corollary 4.12}  In $V (Q | \k \ast Q_\k),$ $E_1$ is
stationary, every stationary $S \subseteq E_0$ reflects fully,
and every
stationary $T \subseteq E_1$ reflects fully.
\endproclaim
 
\demo{Proof}  The first part follows from Lemma 4.3 and
4.10.  The second
part is a consequence of Lemmas 3.8 and 4.6 and the construction
that destroys non-1-positive as well as all non-2-positive sets.
\qed\enddemo
\pbf
 
This concludes Case 2.  We can now go on to Case 3 (and in an
analogous way, to higher cases),  with only one difficulty
remaining.  In analogy with definition 4.2 we can define a filter
$F_3$ and the associated with it 3-positive sets.  All
the proofs of Chapter 4 will generalize from Case 2 to  Case 3,
with the exception of Lemma 4.4 which proved that ``1-positive"
is an absolute $\Pi^1_2$ property.  The
proof does not generalize, as it uses, in an essential way, the
fact that the
property is $\Pi^1_2$, while ``2-positive" is a $\Pi^1_3$
property.
 
However, we can replace the property ``$A \subseteq E_1 \cap \k$
is 2-positive" by another $\Pi^1_3$ property that is
absolute for $\k$, and that is equivalent to the definition 4.2
at all stages of the iteration $Q_\k$ except possibly at the end
of the iteration.  The new
property is as follows:
 
\proclaim{(4.13)}
Either Full Reflection fails for some $S\subseteq Sing \cap \k$
and $A$ is 2-positive, or Full
Reflection holds for all subsets of $Sing$ and $A$ is stationary.
\endproclaim
 
``Full Reflection" for $S \subseteq Sing$ means that $E_0 - Tr
(S)$ is nonstationary.  It is a $\Sigma^1_1$ property of $S$,
and so (4.13) is $\Pi^1_3$.  We claim that Full Reflection
fails at every
intermediate stage of $Q_\k$.  Hence (4.13) is equivalent to
``2-positive" at the intermediate stages.  At the end of $Q_\k$,
every 2-positive set
becomes stationary, and every non-2-positive set becomes
nonstationary.  Hence (4.13) is absolute.
 
Since for every
$\alpha < \k^+$, the size of $Q | \k \ast Q_\k | \alpha$ is $\k$,
the following lemma verifies our claim:
 
\proclaim{Lemma 4.14}  Let $\k$ be a Mahlo cardinal, and assume $V = L
[X]$
where $X \subseteq \k$.  Then there exists a stationary set $S
\subseteq Sing \cap \k$ such that
for every
$\gamma \in E_0,$ $S \cap \gamma$ is nonstationary.
\endproclaim
 
\demo{Proof}  We define $S \subseteq Sing$ by induction on
$\alpha < \k$.
Let $\alpha \in Sing$ and assume $S \cap \alpha$ has been
defined.  Let $\eta (\alpha)$
be the least $\eta < \alpha^+$ such that $L_\eta [X \cap \alpha]$
is a model of ZFC$^-$ and
$L_\eta [X \cap \alpha] \models \alpha$ is not Mahlo.  Let
$$ \alpha \in S \; \text{ iff } \;  L_{\eta (\alpha)} [X \cap
\alpha] \models
S \cap \alpha \; \text{is nonstationary}.$$
First
we show that $S$ is stationary.
 
Assume that $S$ is nonstationary.  Let $\nu < \k^+$ be such that
$S \in L_\nu [X]$ and $L_\nu [X] \models S \text{ is nonstationary.}$
Also, since $\k$ is Mahlo, we have $L_\nu [X] \models \k$ is
Mahlo.
Using a continuous elementary chain of submodels of $L_\nu [X]$,
we find a club $C \subseteq \k$ and a function $\nu (\xi)$ on $C$
such that
for every $\xi \in C$,
$$L_{\nu (\xi)} [X \cap \xi] \models \xi\; \text{is} \;
\text{Mahlo and} \; S \cap \xi \; \text{is not stationary}. $$
If $\alpha \in Sing \cap C$, then
because $\alpha$ is Mahlo in $L_{\nu (\alpha)} [X \cap \alpha]$
but non Mahlo in $L_{\eta (\alpha)} [X \cap \alpha]$, we have
$\nu (\alpha) \le \eta (\alpha)$.  Since $S \cap \alpha$ is
nonstationary in
$L_{\nu (\alpha)} [X \cap \alpha]$, it is nonstationary in
$L_{\eta (\alpha)} [X \cap \alpha]$.
Therefore $\alpha \in S$, and so $S \supseteq Sing \cap C$
contrary to the
assumption that $S$ is nonstationary.
 
Now let $\gamma \in E_0$ be arbitrary and let us show that $S
\cap \gamma$ is nonstationary.
Assume that $S \cap \gamma$ is stationary.  Let $\delta <
\gamma^+$ be such that $S \cap \gamma \in L_\delta [X \cap
\gamma]$, that $L_\delta [X \cap \gamma] \models S \cap \gamma$
is
stationary and that $L_\delta [X \cap \gamma] \models \gamma$ is
not Mahlo.  There is a club $C \subseteq \gamma$ and a function
$\delta (\xi)$ on $C$ such that for every $\xi \in C$,
 
$$ L_{\delta (\xi)} [X \cap \xi] \models \xi \; \text{is not
Mahlo and} \; S \cap \xi \; \text{is stationary}. $$
 
Since $S \cap \gamma$ is stationary, there is an $\alpha \in S
\cap C$.
Because $\eta (\alpha)$ is the least $\eta$ such that $\alpha$ is
not Mahlo in
$L_{\eta (\alpha)} [X \cap \alpha]$, we have $\eta (\alpha) \le
\delta
(\alpha)$.  But $S \cap \alpha$ is nonstationary in $L_{\eta
(\alpha)}
[X \cap \alpha]$ and stationary in $L_{\delta (\alpha)} [X \cap
\alpha]$, a
contradiction.  
\qed\enddemo
 
Now with the modification given by (4.13), the proofs of
Chapter 4 go through in the higher cases, and the proof of
Theorem B is complete.


\Refs
\ref\no 1
\by  J. Baumgartner, L. Harrington and E. Kleinberg
\paper Adding a closed
unbounded set
\jour J.~Symb. Logic 
\vol 41
\yr 1976
\pages 481--482\endref
 
\ref\no2
\by J. Baumgartner, A. Taylor and S. Wagon
\paper  On splitting
stationary subsets of large cardinals
\jour J. Symb. Logic 
\vol 42
\yr 1977
\pages 203--214
\endref
 
\ref\no3 
\by L. Harrington and S. Shelah
\paper  Some exact
equiconsistency results in set theory
\jour Notre Dame J. Formal Logic
\vol 26
\yr 1985 \pages 178--188 \endref
 
\ref\no 4
\by T. Jech
\paper Stationary subsets of inaccessible cardinals
\inbook  Axiomatic Set Theory \ed J. Baumgartner
\bookinfo Contemporary
Math. 
vol. 31
\publ Amer. Math. Soc. \yr 1984 \pages 115--142
\endref
 
\ref\no 5
\by T. Jech
\book Multiple Forcing
\publ Cambridge University
Press
\yr 1986 \endref
 
\ref\no 6
\by T. Jech and H. Woodin
\paper Saturation of the closed
unbounded filter on the set of regular cardinals
\jour Transactions
Amer. Math. Soc. \vol 292 \yr 1985 \pages 345--356
\endref
 
\ref\no 7 \by M. Magidor \paper Reflecting stationary sets
\jour J. Symb. Logic
\vol 47 \yr 1982 \pages 755--771 \endref
\endRefs
\enddocument
\end